\documentclass[11pt]{article}
\usepackage{graphicx}
\usepackage{amssymb, amsmath,}
                                                                                                                                                                                                                                                                                                                                                                                                                                                                                                                                                                                                                                                                                                                                                                                                                                                                                                                                                                                                                                                                                                                                               \newtheorem{theorem}{Theorem}[section]

\newtheorem{definition}[theorem]{Definition}
\newtheorem{lemma}[theorem]{Lemma}
\newtheorem{Remark}[theorem]{Remark}
\newtheorem{Proposition}[theorem]{Proposition}
\title{\textbf{The weight distribution of the functional\\ codes defined by forms of degree 2 on hermitian surfaces}}          

\author{\textbf{Fr\'ed\'eric A. B. Edoukou}  \\ \\ CNRS, Institut de Math\'ematiques de Luminy  \\ Luminy case 907  - 13288  Marseille Cedex 9 - France \\  E.mail : edoukou@iml.univ-mrs.fr }
\begin{document}
\date{}
\maketitle

{\footnotesize \begin{flushleft}
\textbf{Abstract}
\end{flushleft}
We study the functional codes $C_2(X)$ defined on a projective variety $X$, in the case where $X \subset {\mathbb{P}}^3$ is a non-degenerate hermitian surface. We first give some bounds for  $\# X_{ Z(\mathcal{Q})}(\mathbb{F}_{q} )$, which are better than the ones known.
We compute the number of codewords reaching the second weight. We also estimate the third weight, show the geometrical structure of the  codewords reaching this third weight and compute their number.
The paper ends with a conjecture on the fourth weight and the fifth weight of the code $C_2(X)$.
}
\section{Introduction}
The code $C_h(X)$ over $\mathbb{F}_{q}$ ($q=t^2$, with $t$ a prime power) has been studied by S\o rensen  $\lbrack 10 \rbrack$ in his Ph.D Thesis in the case where $X$ is a non-singular hermitian surface and also by G. Lachaud \lbrack 8\rbrack. In fact, the code $C_2(X)$ has previously been studied over $\mathbb{F}_{4}$ by Paul P. Spurr $\lbrack 11 \rbrack$ in his master thesis; by a complete computer search he found the weight distribution of this code.\\
S\o rensen was not able to find the minimum distance, even less the weight distribution. For $f$ a form of degree $h$, he only conjectured the following result: $$\mathrm{for} \ h\le t \quad \# X_{ Z(f)}(\mathbb{F}_{q} ) \le h( t^{3}+ t^{2}-t)+t+1,$$ which should give the minimum distance. In $\lbrack 4\rbrack$ the author didn't only prove the conjecture for $h\le2$, but determined the second weight and three types of quadrics giving codewords of second weight. \\
The purpose of this paper is to express as possible the weight distribution of the code $C_2(X)$. Indeed, the study of the weights is important since they measure the efficiency of the code, and their notification is useful in the procedures of decoding. \\
The paper has been organized as follows. First of all, we give some bounds for the number of points of the section of the hermitian surfaces with quadrics of rank 2, quadrics cone, and the elliptic quadrics. These bounds are better than the ones obtained in  $\lbrack 4 \rbrack$. Thus, codewords of second weight are obtained only from the three types of quadrics found in $\lbrack 4 \rbrack$; we compute their number. We determine the third weight, the geometrical structure of codewords reaching it and we compute their number. Finally, we conclude our work by a conjecture on the fourth and fifth weight of the code $C_2(X)$.\\
We now give a summary of notation and terminology. We denote by $\mathbb{F}_q$ the field with $q$ elements and by $\mathbb{P}^{n}= \mathbb{P}^n(\mathbb{F}_q)=PG(n,q)$ the projective space of $n$ dimension over $\mathbb{F}_q$. We use the term forms of degree $h$ to describe homogeneous polynomials $f$ of degree $h$; and we denote by $\mathcal{Q}=Z(f)$, the zeros of $f$ in $\mathbb{P}^{n}(\mathbb{F}_q)$, which is a quadric when $h=2$.
Let $\mathcal{Q}$ be a quadric, the rank of $\mathcal{Q}$ denoted $r(\mathcal{Q})$ is the smallest number of indeterminates appearing in $f$ under any change of coordinate system.\\
 Let $\mathcal{F}$ be the vector space of forms of degree 2, $X\subset \mathbb{P}^3$ the non-degenerate hermitian surface $X: x_0^{t+1}+x_1^{t+1}+x_2^{t+1}+x_3^{t+1}=0$ and $\vert X\vert$ the number of rational points of $X$ over $\mathbb{F}_q$. We denote by $W_i$ the set of points with homogeneous coordinates $(x_0:...:x_3)\in \mathbb{P}^3$ such that $x_j=0$ for $j<i$ and $x_i=1$. The family $\{W_i\}_{0\le i\le 3}$ is a partition of $\mathbb{P}^3$.
The code $C_2(X)$ is the image of the linear map  
 $c: \mathcal{F}
  \longrightarrow
  \mathbb{F}_{q}^{\vert X\vert}$, defined by $c(\mathcal{Q})={(c_x(\mathcal{Q})})_{x\in X}$,  where  $c_x(\mathcal{Q})= \mathcal{Q}(x_0,...,x_3)$ with $x=(x_0,...,x_3) \in W_i$.
\section{The number of intersection points of the hermitian surface with some quadrics}
The subject of this section is to find some bounds for the number of intersection points of the hermitian surface and some quadrics which are better than the ones of $\lbrack 4\rbrack$. In what follows we will note: 
$$ s(t)=2t^3+2t^2-t+1 \qquad  s_2(t)=2t^3+t^2+1\qquad s_3(t)=2t^3+t^2-t+1$$  
$$s_4(t)=2t^3+1\qquad s_5(t)=2t^3-t+1$$ 
\subsection{Rank ($\mathcal{Q})=2$ and $\mathcal{Q}$ is a pair of distinct non-tangent planes}
If we write $\mathcal {Q } =  \mathcal {P} _{1} \cup   \mathcal {P} _{2}$, then $\mathcal {P }_{1} \cap X $ and $\mathcal {P }_{2} \cap X $ are non-singular hermitian curves in $PG(2,t^2)$, such that $\vert  \mathcal {P }_{1} \cap X \vert  = \vert \mathcal {P }_{2} \cap X \vert =t^3+1$. \\
Let $l= \mathcal {P} _{1} \cap \mathcal {P} _{2}$, then $l$ is not contained in $X$ from theorem 9.1 $\lbrack 2\rbrack$ p.1176. Thus, $l\cap X$ is either a single point or a set of $t+1$ points. Hence $\vert \mathcal {Q } \cap X \vert = s_4(t)$ when $l\cap X$ is a single point or $\vert \mathcal {Q } \cap X \vert  = s_5(t)$ when $l$ is a secant.
\subsection{Rank ($\mathcal{Q}) =3$ ($\mathcal{Q}$ is a cone)}
Here $\mathcal{Q}$ consists of the points on $q+1$ lines passing through a vertex no three of which are coplanar. Thus, the intersection of the cone with the hermitian surface contains at most two generators. If no line of the cone is contained in $X$, we get $\vert \mathcal{Q}\cap X\vert \le t^3+t^2+t+1< s_4(t)$. From proposition 5.3 $\lbrack 4\rbrack$, we conclude that, $\vert \mathcal{Q}\cap X\vert =t^3+t^2+1$ or  $\vert \mathcal{Q}\cap X\vert = t^3+2t^2-t+1$ according to the                                                                         cone contains exactly one or two lines of $X$. 
\subsection{Rank ($\mathcal{Q}) =4$ and $\mathcal{Q}$ is an elliptic quadric }
When $\mathcal {Q }$ is an elliptic quadric, 
no line is contained in $\mathcal {Q }$. 
We get two cases: either $\vert X \cap \mathcal {Q } \vert \le 1$ or $\vert X \cap \mathcal {Q } \vert \ge 2$.\\
In the case  $\vert X \cap \mathcal {Q } \vert  \ge 2$, we will use an analogous technique as the one used in  $\lbrack 4\rbrack$  for the elliptic quadric.\\
Let us choose two distinct points $ P_1, P_2  \in X \cap \mathcal {Q }$.  Let $H$ be a plane passing through the points  $ P_1$ and $P_2$. 
In \lbrack 6\rbrack  \ p.156, table 7.2, J.W.P. Hirschfeld gave the classification of quadrics in $PG(2,q)$. Thus, we get that  $ \mathcal{Q} \cap H$ is either a line or a pair of distinct lines or a point (a pair of conjuguate lines in $PG(2,q)$ which meet in $PG(2,t)$) or a conic (a non-degerate quadric). Since the quadric $ \mathcal{Q}$  does not contain any line and contains two points, we can say that $ \mathcal{Q} \cap H $ is a conic (a non-degerate quadric) and $H$ is non-tangent to $\mathcal{Q}$. \\
Let $(\mathcal{D})$ be the line passing through the two points $P_1$ and $P_2$. 
We consider also  all the planes $(H_{i} )_{i \in I}$ passing through  the line $(\mathcal{D})$. Then  $ \mathcal{Q} \cap H_{i} $ as a non-degenerate quadric is a non-singular curve.
On the other hand $X \cap H_{i} $ is a hermitian curve (singular or not).
Since $H_i$ is non-tangent to $\mathcal{Q}$, as in $\lbrack 4\rbrack$
we conclude that: $$\qquad  \ \mathrm{for \ each} \quad i  \qquad \vert X\cap \mathcal{Q} \cap H_{i}  \vert \le 2(t+1)  \quad (2.3.1)$$
The two points $P_1$ and $P_2$ belong to each  $X \cap \mathcal{Q} \cap H_{i}$, and therefore from (2.3.1) we get: $ \quad \vert X \cap \mathcal{Q} \cap H_{i}  -\{P_1,P_2\}  \vert \le 2(t+1) -2 \quad(2.3.2)$\\
There are exactly $q+1$ planes passing through the line $(\mathcal{D})$, and their union generate the whole projective space $PG(3,q)$. 
Thus, we get:  $$\qquad \qquad \quad \mathcal{Q} \cap X = (X\cap \mathcal{Q}\cap H_1)\cup \bigcup_{i=2}^{q+1}(  X \cap \mathcal{Q} \cap H_{i}  -\{P_1,P_2\} ) \quad (2.3.3)$$ 
And from the relations (2.3.1), (2.3.2) and (2.3.3) we deduce that $ \vert \mathcal{Q} \cap X  \vert \le 2(t+1)+q(2(t+1)-2)$. Thus, we get finally $ \vert \mathcal{Q} \cap X  \vert  \le 2t^3+2t+2$.
We have $ 2t^3+2t+2 < s_2(t)$ for $t\ge 3$ and for $t=2$ we also have $ \vert \mathcal{Q} \cap X  \vert < \vert \mathcal{Q}\vert <s_2(2)$. Then when $\mathcal{Q}$ is an elliptic quadric we conclude that $ \vert \mathcal{Q} \cap X  \vert < s_2(t).$
\section{The second weight of the code $C_{2}(X)$}
\begin{theorem} $\lbrack 4\rbrack$
The second weight of the code $C_{2}(X)$ is  $t^5-t^3$.    
\end{theorem}
\begin{theorem}
The codewords of second weight correspond to:\\ 
\textendash \ hyperbolic quadrics containing three skew lines of $X$,\\
\textendash \  quadrics which are union of two tangent planes meeting at a line contained in $X$,\\
\textendash \ quadrics which are pair of two planes, one tangent to $X$ and the second not tangent to $X$ and the line of intersection of the two planes meeting $X$ at a single point.
\end{theorem}
\textbf{Proof} From paragraph 2.3, we deduce that the codewords of second weight can not correspond to elliptic quadrics. Therefore the three types of quadrics found in $\lbrack 4\rbrack$ are the only ones which give codewords of second weight.  
\begin{lemma}
Let $\mathcal{P}_{1}$ and $\mathcal{P}_{1}^{\prime}$ two tangent planes to $X$, $\mathcal{P}_{2}$ and $\mathcal{P}_{2}^{\prime}$  two non-tangent planes with   $\mathcal{P}_1\cup \mathcal{P}_2 \ne  \mathcal{P}^{\prime}_1\cup \mathcal{P}^{\prime}_2$ and such that $\vert \mathcal{P}_{1}\cap \mathcal{P}_{2}\cap X\vert=\vert \mathcal{P}_{1}^{\prime}\cap  \mathcal{P}_{2}^{\prime}\cap X\vert=1$. Then $(\mathcal{P}_{1}\cap X)\cup (\mathcal{P}_{2}\cap X)\ne (\mathcal{P}_{1}^{\prime} \cap X)\cup   (\mathcal{P}_{2}^{\prime}\cap X)$.
\end{lemma}
\begin{lemma}
Let $\mathcal{Q}^{\prime}$ be a hyperbolic quadric containing three skew lines of $X$. Let $\mathcal{Q}=\mathcal{P}_1\cup \mathcal{P}_2$ be a pair of disctinct planes, one tangent to $X$. Then $\mathcal{Q}^{\prime} \cap X \ne (\mathcal{P}_1\cap X)\cup (\mathcal{P}_2 \cap X)$.
\end{lemma}
\begin {lemma}
Let $\mathcal{Q}=\mathcal{P}_1\cup \mathcal{P}_2$ be a pair of distinct tangent planes and $l=\mathcal{P}_1\cap \mathcal{P}_2 \subset X$. Let $\mathcal{Q}^{\prime}=\mathcal{P}_1^{\prime} \cup \mathcal{P}_2^{\prime}$ be two  planes, $\mathcal{P}_1^{\prime}$ tangent to $X$, $\mathcal{P}_2^{\prime}$ non-tangent to $X$ and $l^{\prime}=\mathcal{P}_1^{\prime}\cap \mathcal{P}_2^{\prime}$ intersecting $X$ at a single point. Then $(\mathcal{P}_1\cap X)\cup (\mathcal{P}_2\cap X)\ne (\mathcal{P}_1^{\prime} \cap X)\cup (\mathcal{P}_2^{\prime} \cap X)$.
\end{lemma}
The proofs of these three lemmas are obvious. There are analogous to the proofs of lemma 6.3 and lemma 6.4 $\lbrack 4\rbrack$.\\\\
We come now to the computation of the codewords of second weight. 
\begin{theorem}
The number of codewords of second weight is: 
$$\frac{1}{2}(t^2-1)(t^5+t^3+t^2+1)(3t^2-t+1)t^2$$
\end{theorem}
\textbf{Proof}
Since there are three types of codewords of second weight  we will find the number of each of them.\\
$\textbf{(i)}$ From J.W.P. Hirschfeld \lbrack 7\rbrack \  p.124, there are exactly $n_{q}=\frac{1}{2}q^2(q\sqrt{q}+1)(q+1)$ hyperbolic quadrics containing three skew lines of the hermitian surface.  Let $\mathcal{Q}$ and $\mathcal{Q}^{\prime}$  be two distinct hyperbolic quadrics giving  codewords of second weight, we have $\mathcal{Q}\cap X \ne \mathcal{Q}^{\prime}\cap X$. Thus, there are exactly $(q-1)\lbrack\frac{1}{2}q^2(q\sqrt{q}+1)(q+1) \rbrack$ codewords of second weight obtained from the hyperbolic quadrics. \\
$\textbf{(ii)}$ We also know from \lbrack 2\rbrack, theorm 7.3 p.1172, that there are exactly $\#X(\mathbb{F}_{q})\\=t^5+t^3+t^2+1$ tangent planes to the hermitian surface X. We want now to determine the number of quadrics $\mathcal{Q}=\mathcal{P}_{1}\cup \mathcal{P}_{2}$ where $\mathcal{P}_{1}$ and $\mathcal{P}_{2}$ are two distinct tangent planes meeting at a line contained in $X$. Given a tangent plane $\mathcal{P}_{1}$, there are exactly $t^2(t+1)=t^3+t^2$ possibilities to choose the plane $\mathcal{P}_{2}$; so we get $(t^5+t^3+t^2+1)(t^3+t^2)$ couples $(\mathcal{P}_{1},\  \mathcal{P}_{2})$. Since the two couples $(\mathcal{P}_{1},\  \mathcal{P}_{2})$ and $(\mathcal{P}_{2},\  \mathcal{P}_{1})$ are the same (they give the same quadric), we get from these tangent planes $\frac{1}{2}(t^5+t^3+t^2+1)(t^3+t^2)$ quadrics giving codewords of second weight. And from lemma 6.3 and 6.4 $\lbrack 4\rbrack$, the codewords obtained from these quadrics are all distinct. Thus, we get $(t^2-1)\lbrack \frac{1}{2}(t^5+t^3+t^2+1)(t^3+t^2)\rbrack$ codewords of second weight.\\
$\textbf{(iii)}$ Let $Q=\mathcal{P}_{1}\cup \mathcal{P}_{2}$ be a quadric where $\mathcal{P}_{1}$ and $\mathcal{P}_{2}$ are respectively a tangent plane to $X$, a non-tangent plane to $X$, where $l=\mathcal{P}_{1}\cap \mathcal{P}_{2}$ with $l\cap X=\{P_1\}$ ($l \not \subset X$).
Here, we need to find first of all the number of planes $(H_i)_{i\in I}$ for which  we have 
$H_i \cap \mathcal{P}_1=l$ where $l$ is a line intersecting X in a single point $P_1$.\\
We know that there are $(t^2-t)$ lines of such type (lines $l$).\\
Given a line $l$, there pass exactly $q+1$ planes through $l$  ($\mathcal{P}_1$ is one of these planes). Thus, we get at maximum $q(t^2-t)$ planes $H_i$,  i.e. $t^4-t^3$ planes $H_i$. In fact, there are exactly $t^4-t^3$ planes $H_i$, since if $H_1$ and $H_2$ are respectively planes through $l_1$ and $l_2$ two distinct lines, then $H_1\ne H_2$.\\
We also know that the number of tangent planes to $X$ is $\#X((\mathbb{F}_{q})=t^5+t^3+t^2+1$, so we deduce  that there are exactly $(t^5+t^3+t^2+1)(t^4-t^3)$ quadrics $Q=\mathcal{P}_{1}\cup \mathcal{P}_{2}$ (where $\mathcal{P}_{1}$ is a tangent plane to $X$, $\mathcal{P}_{2}$ is a non-tangent plane to $X$, and $l=\mathcal{P}_{1}\cap \mathcal{P}_{2}$ intersecting $X$ at a single point). And from lemma 3.3, the codewords obtained from these quadrics are all distinct. Thus, we get exactly $(t^2-1)(t^5+t^3+t^2+1)(t^4-t^3)$ codewords of second weight from these quadrics.\\  
 From lemmas 3.4 and 3.5, we deduce that the codewords obtained from $\textbf{(i)}$, $\textbf{(ii)}$ and $\textbf{(iii)}$ are all distinct.\\ Therefore we have exactly:  $(t^2-1)\lbrack\frac{1}{2}t^4(t^3+1)(t^2+1)\rbrack+(t^2-1)\lbrack\frac{1}{2}(t^5+t^3+t^2+1)(t^3+t^2)\rbrack+(t^2-1)(t^5+t^3+t^2+1)(t^4-t^3)$ $=\frac{1}{2}(t^2-1)(t^5+t^3+t^2+1)(3t^2-t+1)t^2$ codewords of second weight.\\\\
 Observe that, in the case $t=2$, there are exactly $\frac{1}{2}\times 3\times 45\times 11\times 4=2970$ codewords of second weight. Thus, we recover the result of Paul P. Spurr \lbrack 3\rbrack  \ p.120 calculated by computer.
\section{The Weight distribution of the code $C_2(X)$}
In what follows 
$s(t)$, $s_2(t)$, $s_3(t)$, $s_4(t)$ and $s_5(t)$ have the same values as in section 2. It has been shown in $\lbrack 4\rbrack$ with the help of the relations (3.2.2) and (3.4.3) that $s(t)$ and $s_2(t)$ give respectively the first weight (minimum distance) and the second weight. \\
For $t > 3$,  $s_3(t) > 2t^3+2t+2$; and for $t=2$, the elliptic quadric $\mathcal{Q}$ over $\mathbb{F}_{4}$ has $\vert \mathcal{Q}\vert = 17$ points, and  we have $17< s_3(2)=19$. Therefore, from table 4.2, the third weight is given by $s_3(t)$ for $t\ne 3$. 
\begin{theorem}
Let $l$ be a line in $PG(3, \mathbb{F}_{q})$ and $X$ the hermitian surface
then  $l$ meets $X$ at a single point or $t+1$ points or is contained in $X$. This line is respectively called a tangent, a secant, or a generator. 
\end{theorem}
\textbf{Proof:}  See $\lbrack 2\rbrack$, p.1179.
\begin{definition}
 A regulus (notion used by Hirschfeld \lbrack 7\rbrack \  p.4) is the set of transversals of three skew lines. It consists of $q+1$ skew lines. 
\end{definition}
\begin{theorem} \lbrack 7\rbrack \ p.23
A hyperbolic quadric $\mathcal{Q}$ consists of $(q+1)^2$ points, which are all on a pair of complementary reguli. The two reguli are the two systems of generators of $\mathcal{Q}$.  
A hyperbolic quadric whose complementary reguli are $\mathcal{R}$ and  $\mathcal{R}^{\prime} $ is denoted by  $\mathcal{H}( \mathcal{R} , \mathcal{R}^{\prime} )$.
\end{theorem}
\begin{Remark}
When $\mathcal{Q}$ is a hyperbolic quadric containing exactly two skew lines on the hermitian surface $X$ (i.e. type 12 in table 4.2), we  get: \\
\textendash\ for $t\ge 3$, $\# X_{Z(\mathcal{Q})}(\mathbb{F}_{q}) \le  t^3+3t^2-t+1 < s_4(t)$,\\
\textendash\ for $t=2$,  we get $\# X_{Z(\mathcal{Q})}(\mathbb{F}_{4}) \le  2^3+3\times 2^2-2+1 =19= s_3(2)$. Let us suppose in the case $(t=2)$ that, $\# X_{Z(\mathcal{Q})}(\mathbb{F}_{4}) < 19$;  therefore it exists a line $l$ on the considered regulus of $\mathcal{Q}$ such that $l\cap X$ is a single point. Thus we get $\# X_{Z(\mathcal{Q})}(\mathbb{F}_{4}) \le  2(t^2+1)+(t^2-2)(t+1)+1 =17 = s_4(2)$. Therefore we can assert that $$\# X_{Z(\mathcal{Q})}(\mathbb{F}_{4}) = s_3(2)\ \  \mathrm{or}\ \  \# X_{Z(\mathcal{Q})}(\mathbb{F}_{q}) \le  s_4(t) \qquad(4.4.1)$$
\end{Remark}
\begin{theorem}
Let $\mathcal{Q}$ be a hyperbolic quadric containing exactly two skew lines on the hermitian surface $X$. Then we have $\# X_{Z(\mathcal{Q})}(\mathbb{F}_{q}) \le  s_4(t)$.
\end{theorem}
\textbf{Proof }  see appendix.
\subsection{Table of the intersection of quadrics with the hermitian surface}
Let us show in a table the geometrical structure of the intersection of quadrics with the hermitian surface. \\

\hspace{15mm}
\begin{tabular}{|c|c|c|}

	\hline
	Rank ($\mathcal{Q}$) & Description & Type  \\
  	\hline
	\hline
	1 		      & repeated plane 			  &$1, 2 $				 							       						  \\
	\hline
	2  & line	  	  & $ 3,4,5$					 							       					\\
	
	\hline		
	2 		      & pair of distinct planes 			  & $ 6,7,8$				 							       		 \\
	\hline
	3	      & quadric cone 			  & $ 9,10$				 							       			\\
	\hline
	 	4	       & hyperbolic quadric 			  & $11, 12, 13, 14$				 							       				 \\
	\hline

	 4             &elliptic quadric                               &$15$                                                         
	
	                                    \\
	
       \hline            
	\end{tabular}

\vspace{3mm}

\begin{tabular}{cc}
\hspace{-14mm}

\begin{tabular}{|c|c|}

	\hline
	Type & Description    \\
  	\hline
	\hline
	1 		      &  a tangent plane to $X$			 			 							       						  \\
	\hline
	2  & a non-tangent plane to $X$	  	 					 							       					\\
	
	\hline		
	3 		      & tangent to $X$  			  				 							       		 \\
	\hline
	4	      & secant  			  			 							       			\\
	\hline
	5	       & a generator of $X$ 						 							       				 \\
	\hline

	 6             &two tangent planes to $X$                                                                                
	
	                                    \\
	
       \hline  
       
     7             &$\mathcal{P}_1$ is non-tangent to $X$  \\                                                                     
        
        &and $\mathcal{P}_2$ is tangent to $X$ 
	
	                                    \\
	
       \hline  
            
	\end{tabular}
&
\begin{tabular}{|c|c|}

	\hline
	Type & Description    \\
  	\hline
	\hline
	8 & two non-tangent planes    \\
  	
	\hline
	9 		      &no generartor in the cone 			 			 							       						  \\
	\hline
	10  & 1 or 2 generators in the cone  	 					 							       					\\
	
	\hline		
	11 		      & 3 skew generators in $\mathcal{R}$			  				 							       		 \\
	\hline
	12	      & 2 skew generators  $\mathcal{R}$  			  			 							       			\\
	\hline
	13	       & 1 generator in $\mathcal{R}$ 					 							       				 \\
	\hline

	 14             & no generator in $\mathcal{R}$                                                                              
	
	                                    \\
	
       \hline

	        15            &elliptic quadric                                                                                
	
	                                    \\
	
       \hline  
             
	\end{tabular}

\end{tabular}
\vspace{3mm}

\subsection{Table of the weight distribution}
Let us show in a table the quadrics (geometrical structure) with the corresponding weights 
associated to the codewords obtained from them. From the results of section 2  and those of section 5 $\lbrack 4\rbrack$, we have the following table.\\

\hspace{-21mm}
\begin{tabular}{|c|c|l|c|l|c|} 
	\hline
	Rank & Type& $\# X_{Z(\mathcal{Q})}(\mathbb{F}_{q})$ & $\# X_{Z(\mathcal{Q})}(\mathbb{F}_{4})$ & Weight over $\mathbb{F}_{q}$ & Weight over $\mathbb{F}_{4}$ \\
  	\hline
	\hline
	1 		      & 1 			  &$ t^3+t^2+1 $				& 13 							       &$ t^5	$						 & 32 \\
	\cline{2-6}
	(repeated plane) & 2 	  	  & $ t^3+1  $					& 9 							       & $ t^5+t^2 $					& 36 \\
	\hline		
	2 		      & 3 			  & $ 1 $					& 1 							       &$  t^5+t^3+t^2 $				& 44 \\
	\cline{2-6}
	(line)		      & 4 			  & $ t+1  $					& 3 							       & $ t^5+t^3+t^2-t $				& 42 \\
	\cline{2-6}
	 		       & 5 			  & $ t^2+1 $					& 5 							       &$  t^5+t^3 $					& 40 \\
	\hline

	 {}              &6                               &$s_4(t)$                                                      & 17      
	
	                      &$t^5-t^3+t^2$                & 28\\
	
	\cline{3-6} 
	2                 &{ }                         & $s_5(t) $                                       &15                  &$t^5-t^3+t^2+t$             &  30  \\   
	
	\cline {2-6}  
	(pair of dis-    & {}                    & $s_3(t)$                                                                & 19
	                   &$t^5-t^3+t$          &  36     \\
	\cline {3-6}
	tinct planes)     &7                       & $s_2(t)$                                                                & 21
	                 &$t^5-t^3$            & 24          \\

	\cline{2-6} 
	$\mathcal{Q}=\mathcal{P}_1\cup \mathcal{P}_2$             &{}                        &$s(t)$                                                                      &23                                                                               
	                &$t^5-t^3-t^2+t$     & 22   \\
	 \cline{3-6} 
	  {}      &  8                             &$s_2(t)$                                                                 &21
	           &$t^5-t^3$                 &24   \\       
   \hline            
	3         & 9                             & $\le t^3+t^2+t$                                  & $\le 15$    
	           & $\ge t^5-t$           & $ 30 \le w \le 44$  \\   
	{}        & {}                              & $+1<  s_4(t)$                                  &{}
	           & {}                             & {} \\
  \cline{2-6}
      (cone)    & 10                          &$t^3+t^2+1$                                                             &$13$     
                      &$t^5$    &$ 32 $  \\
   
   \cline{3-6}
      {}             & {}                            &$t^3+2t^2-t+1$                                                      &$15$
                      &$ t^5-t^2+t$     &$30$\\

    \hline     
    \hline
        {}          & 11                          &$s_2(t)$                                                                     &21
                    &$t^5-t^3$          & 24	   \\          
	           
   \cline{2-6}
       4       &   12                          &$\le t^3+3t^2-t$                                 &$\le 19$                                                                             &$\ge t^5-t^3$                              & $26\le w\le 32$ \\
        {}   &    {}                                &  $+1 \le  s_3(t)$                                                      & {}
        &  {}                                         & {} \\                  
      \cline{2-6}                              
       (hyperbolic) &    13              &  $\le t^3+2t^2$                                & $\le 17$        &$\ge t^5-t^2$
                              & $28\le w \le 36$      \\   
          {$\mathcal{H}( \mathcal{R}, \mathcal{R}^{\prime})$}                  &     {}               & $+1  \le  s_4(t) $                            & {}                      & {}
                                & {} \\                   	          
       \cline{2-6}           
	           {}      &     14              & $\le t^3+t^2+$                                 & $\le 15$                          & $\ge t^5-t$                                  & $30\le w \le 40 $ \\
	           {}      &      {}                & $ t+1 < s_4(t)$                                                  & {}                              
	           & {}        & {} \\
    \hline 
     4 &         {}             & $\le 2t^3+2t$                &  {}                                                                        &$\ge t^5-t^3+t^2$                                &{}\\

   (elliptic)               & 15               &  $+2 <  s_2(t)$                                 & $\le 17$                  & $- 2t-1$     & $28 \le w \le 45 $ \\
     
    \hline
\end{tabular}

\vspace{4mm}
Note that, for the code $C_2(X)$ over $\mathbb{F}_{4}$, the codewords corresponding to quadrics of type.1 to type.14 have even weights. It has also been proved by Spurr (by a computer program) that  $C_2(X)$ over $\mathbb{F}_{4}$ is an even weight code.
\section{The third weight of the code $C_2(X)$}
\begin{theorem}
The third weight (for $t\ne 3$)of the code $C_2(X)$ is $t^5-t^3+t$
 \end{theorem}
 \textbf{Proof}
 In fact when $\mathcal{Q}$ is a quadric not corresponding to the codewords of minimum weight and second weight, we get $\# X_{Z(\mathcal{Q})}(\mathbb{F}_{q})\le s_3(t)$.\\\\
Note that in the case $t=2$, we recover the third weight of the code $C_2(X)$ over $\mathbb{F}_{4}$ in agreement with \lbrack 3\rbrack  \  which is $2^5-2^3+2=26$.
 \begin{theorem}
 The codewords of third weight (for $t\ne 3$) correspond to quadrics which are union of two planes, one tangent to the hermitian surface $X$, the second not tangent to $X$, and the line of intersection of the two planes intersects the hermitian surface in $t+1$ points.
 \end{theorem}
 \textbf{Proof}
 The codewords of third weight are not given by hyperbolic quadrics of type 12 in table 4.2 from theorem 4.5, nor by elliptic quadrics if $t\ne 3$.\\\\
We come now to the computation of the codewords of third weight.
 \begin{theorem}
 The number of codewords of third weight (for $t\ne 3$) is: $$(t^2-1)(t^5+t^3+t^2+1)(t^6-t^5)$$
  \end{theorem}
 \textbf{Proof}
 Let us write $\mathcal{Q}=\mathcal{P}_1\cup \mathcal{P}_2$, where:\\
\textendash \ $\mathcal{P}_1$ is a tangent plane to the hermitian surface $X$ at a point $P_1$,\\ 
\textendash \ $\mathcal{P}_2$ is a non-tangent plane to the hermitian surface $X$ ,\\
Ñand $l= \mathcal{P}_1 \cap \mathcal{P}_2$ meets the hermitian surface in $t+1$ points.\\
 Given a fixed point $P_1$ in $\mathcal{P}_1$, there are exactly $q+1$ lines passing through $P_1$ and contained in the same plane. We also know that in a plane, there pass exactly $q^2+q+1$ lines. We wish that the line $l$ of intersection is not a generator, nor a tangent (the number of generators and tangent lines is $q+1$). Thus, the remaining lines which number is $q^2+q+1-(q+1)=q^2$, is the number of possibilities we can choose the line $l$ such that $l$ intersects $X$ at $t+1$ points. 
For each line $l$, it pass exactly $q+1$ planes $(\mathcal{P}_i)_{1\le i \le q+1}$ (counting $\mathcal{P}_1$). Let us consider the $q$ planes $(\mathcal{P}_i)_{2\le i \le q+1}$; they are not all non-tangent to $X$ (some of them are tangent to $X$). For a fixed point $P_1$, there are $q^2.q=q^3$ planes  $\mathcal{P}_i$ constructed and these $q^3$ planes are all distinct.\\
In total there are exactly $(t^5+t^3+t^2+1)q^3$ couples $(\mathcal{P}_1, \mathcal{P}_2)$, where
$\mathcal{P}_1$ is a tangent plane to the hermitian surface $X$, $\mathcal{P}_2$ is a non-tangent plane to the hermitian surface $X$ with $l= \mathcal{P}_1 \cap \mathcal{P}_2$ intersecting the hermitian surface at $t+1$ points. We also know from the proof of theorem 6.5 $\lbrack 4\rbrack$ that $(t^5+t^3+t^2+1)t^5$ is the number of couples $(\mathcal{P}_1, \mathcal{P}_2)$ with $\mathcal{P}_1$ and $\mathcal{P}_2$ two tangent planes to the hermitian surface $X$, with the intersection line meeting $X$ at $t+1$ points. Thus, one deduces that $(t^5+t^3+t^2+1)(t^6-t^5)$ is the number of couples $(\mathcal{P}_1, \mathcal{P}_2)$ giving codewords of third weight. Therefore, there are $(t^2-1)(t^5+t^3+t^2+1)(t^6-t^5)$ codewords of third weight.\\\\
Observe that, for $t=2$, there are exactly $3\times 45 \times 32=4320$ codewords of third weight in agreement with Spurr \lbrack 3\rbrack \ p.120.
\begin{Remark} 
\textbf{Conjecture on the fourth and fifth weight.}\\
The author has tried to study again $\# X_{Z(\mathcal{Q})}(\mathbb{F}_{q})$ for $\mathcal{Q}$ an elliptic quadric and conjecture that $s_4(t)$, $s_5(t)$ should give respectively the fourth and the fifth weight. 
The codewords of fourth weight correspond now (for $t\ge 3$) to:\\
\textendash \ quadrics which are union of two non-tangent planes to $X$, and the line of intersection of the two planes meets $X$ at a single point.\\
The codewords of fifth weight correspond (for $t>3$) to:\\
\textendash \ quadrics which are union of two non-tangent planes to $X$, and the line of intersection of the two planes meets $X$ in $t+1$ points,\\
\textendash \ and some particular elliptic quadrics.\\
Unfortunately no proof is found yet. 
\end{Remark}
\section{Appendix}
\begin{lemma}
Let $\overline{l}$ be a line defined over $\overline{\mathbb{F}}_{q}$ in $PG(3, \overline{\mathbb{F}}_{q})$
and $l=\overline{l}\cap PG(3, \mathbb{F}_{q})$. Then $l$ is empty, or a single point or a line defined over ${\mathbb{F}_{q}}$.
\end{lemma}
\textbf{Proof} It is a direct consequence of the equation of a line $PG(3,q)$.\\\\
Let $\mathcal{Q}$ be a non-degenerate quadric over $\overline{\mathbb{F}}_{q}$. Note $d_i=\ $degree of the curve $\mathcal{C}_{i}$ defined on $\mathcal{Q}$. From $ \lbrack 9\rbrack$ (chap.4  example 2 p.241), we deduce an interesting relation between the arithmetical genus of a curve and its degree. Thus, we get the following result. 
\begin{Proposition}
Let $\mathcal{Q} \subset {\mathbb{P}}^3$ be a non-singular quadric surface in ${\mathbb{P}}^3$ and
 $\mathcal{C}$ an irreducible curve on $\mathcal{Q}$.\\
\textendash If deg $\mathcal{C}=2l+2$ then, $p_a(\mathcal{C})\le l^2\quad (\star)$\\
\textendash If deg $\mathcal{C}=2l+1$ then, $p_a(\mathcal{C})\le l(l-1)\quad (\star \star)$
\end{Proposition}
Let us recall the result of Yves Aubry and Marc Perret on the Weil theorem for singular curves, see $ \lbrack 1\rbrack$. With this result we get a bound for the number of rational points of any irreducible curve not necessary smooth.
\begin{Proposition}
Let $C$ be a reduced connected projective algebraic curve over a finite field $\mathbb{F}_{q}$, with $r$ irreducible components and of arithmetical genus $p_a(C)$.
Then   $ \vert \#C(\mathbb{F}_{q})-(rq+1) \vert \le 2p_a(C)\sqrt{q}$.
\end{Proposition}
\textbf{Proof of theorem 4.5}
From the relation $4.4.1$ of remark 4.4, we need to prove that $\# X_{Z(\mathcal{Q})}(\mathbb{F}_{4}) \ne  s_3(2)$.  A hyperbolic quadric over $\mathbb{F}_{4}$ is a pair of two reguli, each regulus containing five skew lines.
Let $\mathcal{R}$ be the regulus of $\mathcal{Q}$ containing the two skew lines $l_1$ and $l_2$ of the hermitian surface.
 
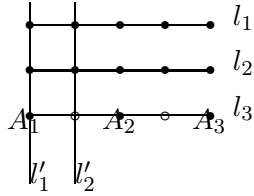
\unitlength=0.6cm
\begin{picture}(7, 5)
\put(1,4){\line(1,0){4}}   \put(1,4){\circle*{.15}}  \put(2,4){\circle*{.15}}  \put(3,4){\circle*{.15}}  \put(4,4){\circle*{.15}}  \put(5,4){\circle*{.15}}  \put(5.5,4){$l_1$}

\put(1,3){\line(1,0){4}}  \put(1,3){\circle*{.15}}  \put(2,3){\circle*{.15}}  \put(3,3){\circle*{.15}}  \put(4,3){\circle*{.15}}  \put(5,3){\circle*{.15}}  \put(5.5,3){$l_2$}

\put(1,2){\line(1,0){4}}  \put(1,2){\circle*{.15}}  \put(2,2){\circle{.15}}  \put(3,2){\circle*{.15}} 
 \put(4,2){\circle{.15}}  \put(5,2){\circle*{.15}}  \put(5.5,2){$l_3$}
 \put(0.5,1.7){$A_1$}  \put(2.6,1.7){$A_2$} \put(4.6,1.7){$A_3$}

\put(1,0.5){\line(0,1){4}} \put(1,0.5){$l_1^{\prime}$}

\put(2,0.5){\line(0,1){4}}  \put(2,.5){$l_2^{\prime}$}

\end{picture}

If each one of the remaining three lines of $\mathcal{R}$ is not  secant to $X$ (i.e. meets $X$ at three points), then the problem is solved. Suppose now that there is a secant line. Let $l_3$ be a line of $\mathcal{R}$ secant to $X$;
 $l_3$ meets $X$ at three points $A_1$, $A_2$ and $A_3$. Let us consider the five planes passing through $l_3$; each of them is tangent to the quadric $\mathcal{Q}$ (since $l_3\cap X=\{A_1, A_2, A_3\}\subset H\cap \mathcal{Q}$). In fact, $H\cap \mathcal{Q}$ as a degenerate quadric of rank 2 in $PG(2,q)$, is a set of 2 lines from J.W.P.Hirschfeld $\lbrack6\rbrack$ p.156. Thus, we can distinguish two cases: $H\cap \mathcal{Q}=\{l_3, l\}$  where $l$ is a line of $\mathcal{Q}$ passing through one of the three points of  $\{A_1, A_2, A_3\}$ or $l$ does not meet any of the points $A_1$, $A_2$ and $A_3$. Without loss of generality, we may assume that                                                                                                                                      $H\cap \mathcal{Q}=\{l_3, l_1^{\prime}\}$ or  $H\cap \mathcal{Q}=\{l_3, l_2^{\prime}\}$. \\\\
\textbf{First case: Suppose that $H \cap \mathcal{Q}=\{l_3, l_1^{\prime}\}$}\\
 If  \textbf{}$H$ is tangent to $X$:{}
In this case $H\cap X$ is a degenerate hermitian curve of rank 2 consisting of 13 points (a set of three lines meeting at a  common point). But, two of the three points $A_1$, $A_2$ and $A_3$ can not be on a the same line of $H\cap X$, otherwise $l_3$ should be contained in $X$. So we get two eventualities:

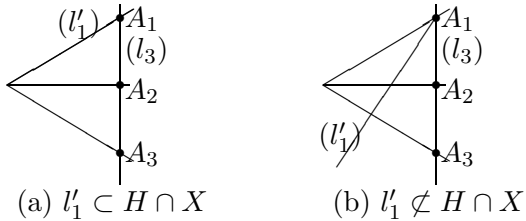
\unitlength=0.6cm
\begin{picture}(12,5)

\put(0.5,3){\thinlines \line(1,0){2.7}} \put(3,3){\circle*{0.15}}  \put(3.1,2.7){$A_2$}
  
  \put(3,0.8){\line(0,1){4}} \put(3.1,3.6){$(l_3)$}
  
  \put(0.5,3){\line(5,3){2.8}}   \put(3,4.5){\circle*{0.15}}  \put(3.1,4.3){$A_1$}
  
  \put(0.5,3){\line(5,-3){2.8}} \put(3,1.5){\circle*{0.15}}  \put(3.1,1.3){$A_3$}
  
  \put(0.5,3){\line(5,3){2.8}}   \put(1.6,4.2){($l_1^{\prime}$)}

  \put(0.7,0.2){(a) $l_1^{\prime} \subset H\cap X$}

 \put(7.5,3){\thinlines \line(1,0){2.7}} \put(10,3){\circle*{0.15}}  \put(10.1,2.7){$A_2$}
  
  \put(10,0.8){\line(0,1){4}} \put(10.1,3.6){$(l_3)$}
  
  \put(7.5,3){\line(5,3){2.8}}   \put(10,4.5){\circle*{0.15}}  \put(10.1,4.3){$A_1$}
  
 \put(7.5,3){\line(5,-3){2.8}} \put(10,1.5){\circle*{0.15}}  \put(10.1,1.3){$A_3$}

 \put(10,4.5){\line(-2,-3){2.2}}  \put(7.4,1.7){($l_1^{\prime}$)}

\put(7.7,0.2){(b) $l_1^{\prime} \not \subset H\cap X$}
\end{picture}

(a) $ l_1^{\prime} \subset H\cap X$:
In this case, $\mathcal{Q}$ contains three lines of $X$ among which two lines are in the same regulus and the third is in the complementary regulus. 
We can now write: $\mathcal{C}=\mathcal{Q}\cap X=l_1\cup l_2\cup l_1^{\prime} \cup \overline{C}$ 
(where $\overline{C}$ is a curve of degree 3 over $\overline{\mathbb{F}}_{4}$). If $\overline{C}$ is irreducible, it is defined on $\mathbb{F}_{4}$ and from proposition 6.2, we have $p_a(\overline{C})=0$ and therefore from proposition 6.3, we get  $\vert \mathcal{Q}\cap X \vert \le 18$. If $\overline{C}$ is not irreducible and since there is no plane conic on $X$, we deduce that $\overline{C}= \overline{C}_1\cup\overline{C}_2\cup \overline{C}_2$ (union of three lines); at most one line among them should be defined over $\mathbb{F}_q$.  Thus, from lemma 6.1 we deduce that $\vert \mathcal{Q}\cap X \vert \le 18$.\\
 (b) $l_1^{\prime} \not \subset H \cap X$:
Here $l_1^{\prime} \cap H \cap X = 3$ points. So that  $\vert H\cap X \cap \mathcal{Q} \vert \le 5$.\\
\textbf{} If H is non-tangent to $X$:{}
In this case, $H\cap X$ is a non-degenerate hermitian curve with $\vert H\cap X\vert = 9$.  Here we need to distinguish two cases:\\
\textendash \ If  $l_1^{\prime} \subset H \cap X $, we recover the case which has been done previously. \\ \textendash \ If $ l_1^{\prime} \not \subset H \cap X$, then $\vert l\cap(H\cap X)\vert = 3$; so  $\vert \mathcal{Q}\cap H \cap X\vert \le 3+3= 6$.\\\\
\textbf{Second case:  Suppose that $H\cap \mathcal{Q}=\{l_3, l^{\prime}_{2}\}$}\\
 Similary, if $H$ is tangent to $X$, then $H\cap X$ is a set of three lines meeting at a common point,  thus $\vert \mathcal{Q}\cap H \cap X\vert \le 6$.
If $H$ is non-tangent to $X$, then $H\cap X$ is a non-singular plane curve and therefore irreducible. Here $H\cap X$ does not contain the line $l_2^{\prime}$, so $l_2^{\prime}$ meets $H\cap X$ in exactly three points. Then $\vert \mathcal{Q}\cap H \cap X\vert \le  6$.\\\\
\textbf{Summary:}
{\itshape When  $\mathcal{Q}$ is a hyperbolic quadric which regulus $\mathcal{R}$ contains exactly two skew lines of the hermitian surface $X$,  $\mathcal{Q}$ should contains an other line $l$ of $X$, in this case $l$ is contained in the complementary regulus of $\mathcal{R}$ and then $\vert \mathcal{Q}\cap X \vert \le18$; or  $\vert \mathcal{Q}\cap H_i\cap X\vert \le 6$ for each one of the five planes $(H_i)_{1\le i \le 5}$ passing through the line $l_3$, which leads to $\vert \mathcal{Q}\cap X \vert \le 18$, therefore $\vert \mathcal{Q} \cap X \vert \in \{13, 15, 17\}$. This concludes the proof of theorem 4.5.}\\\\
 \textbf{Acknowledgment:}
The author would like to thank Mr F. Rodier whose remarks and patience encouraged him to work on the problem. \\\\
\textbf{References}\\
{\footnotesize \lbrack 1\rbrack   \ Y. Aubry and M. Perret, Connected projective algebraic curves over finite fields. Preprint, IML, Marseille, France, 2002.\\\
\lbrack 2\rbrack  \ R. C. Bose and I. M Chakravarti, Hermitian varieties in finite projective space $PG(N,q)$. Canadian J. of Math.18 (1966), pp1161-1182.\\
\lbrack 3\rbrack  \ I. M. Chakravarti, The generalized Goppa codes and related discrete designs from hermitian surfaces in $PG(3, s^2)$. Lecture Notes in computer Sci 311. (1986), pp 116-124.\\
\lbrack 4\rbrack  \ F. A. B. Edoukou, Codes defined by forms of degree 2 on hermitian surface and S\o rensen conjecture. 16 pp., Submitted to Finite Fields and Their Applications (2005).\\
\lbrack 5\rbrack \ R. Hartshorne, Algebraic geometry, Graduate texts in mathematics 52, Springer-Verlag, 1977. \\
\lbrack 6\rbrack  \ J. W. P. Hirschfeld, Projective Geometries Over Finite Fields (Second Edition) Clarendon  Press. Oxford 1998.\\
\lbrack 7\rbrack  \ J. W. P. Hirschfeld, Finite projective spaces of three dimensions, Clarendon press. Oxford 1985. \\
\lbrack 8\rbrack  \ G. Lachaud, Number of points of plane sections and linear codes defined on algebraic varieties;  in "Arithmetic, Geometry, and Coding Theory". (Luminy, France, 1993), Walter de Gruyter, Berlin-New York, 1996, pp 77-104. \\
\lbrack 9\rbrack \  I. R. Shafarevich, Basic algebraic geometry 1, Springer-Verlag, 1994.\\
\lbrack 10\rbrack  \ A. B. S\o rensen, Rational points on hypersurfaces, Reed-Muller codes and algebraic-geometric codes. Ph. D. Thesis, Aarhus, Denmark, 1991.\\
\lbrack 11\rbrack  \ P. P. Spurr, Linear codes over $GF(4)$. Master's Thesis, University of North Carolina at Chapell Hill, USA, 1986.}
\end{document}